%
%
%
\def\unredoffs{} \def\redoffs{\voffset=-.31truein\hoffset=-.59truein}
\def\speclscape{\special{ps: landscape}}
%
%
%
%
\newbox\leftpage \newdimen\fullhsize \newdimen\hstitle \newdimen\hsbody
\tolerance=1000\hfuzz=2pt
\catcode`\@=11 
\def\bigans{b }
\def\answ{b }
\ifx\answ\bigans\message{(This will come out unreduced.}
\magnification=1200\unredoffs\baselineskip=16pt plus 2pt minus 1pt
\hsbody=\hsize \hstitle=\hsize 
\else\message{(This will be reduced.} \let\l@r=L
\magnification=1000\baselineskip=16pt plus 2pt minus 1pt \vsize=7truein
\redoffs \hstitle=8truein\hsbody=4.75truein\fullhsize=10truein\hsize=\hsbody
\output={\ifnum\pageno=0 
  \shipout\vbox{\speclscape{\hsize\fullhsize\makeheadline}
    \hbox to \fullhsize{\hfill\pagebody\hfill}}\advancepageno
  \else
  \almostshipout{\leftline{\vbox{\pagebody\makefootline}}}\advancepageno 
  \fi}
\def\almostshipout#1{\if L\l@r \count1=1 \message{[\the\count0.\the\count1]}
      \global\setbox\leftpage=#1 \global\let\l@r=R
 \else \count1=2
  \shipout\vbox{\speclscape{\hsize\fullhsize\makeheadline}
      \hbox to\fullhsize{\box\leftpage\hfil#1}}  \global\let\l@r=L\fi}
\fi
%
\newcount\yearltd\yearltd=\year\advance\yearltd by -1900

\def\Title#1#2{\nopagenumbers\abstractfont\hsize=\hstitle\rightline{#1}%
\vskip 1in\centerline{\titlefont #2}\abstractfont\vskip .5in\pageno=0}
\def\Date#1{\vfill\leftline{#1}\tenpoint\supereject\global\hsize=\hsbody%
\footline={\hss\tenrm\folio\hss}}
%

\def\draftmode{\message{ DRAFTMODE }\def\draftdate{{\rm preliminary draft:
\number\month/\number\day/\number\yearltd\ \ \hourmin}}%
\headline={\hfil\draftdate}\writelabels\baselineskip=20pt plus 2pt minus 2pt
 {\count255=\time\divide\count255 by 60 \xdef\hourmin{\number\count255}
  \multiply\count255 by-60\advance\count255 by\time
  \xdef\hourmin{\hourmin:\ifnum\count255<10 0\fi\the\count255}}}
\def\nolabels{\def\wrlabeL##1{}\def\eqlabeL##1{}\def\reflabeL##1{}}
\def\writelabels{\def\wrlabeL##1{\leavevmode\vadjust{\rlap{\smash%
{\line{{\escapechar=` \hfill\rlap{\sevenrm\hskip.03in\string##1}}}}}}}%
\def\eqlabeL##1{{\escapechar-1\rlap{\sevenrm\hskip.05in\string##1}}}%
\def\reflabeL##1{\noexpand\llap{\noexpand\sevenrm\string\string\string##1}}}
\nolabels
%
\global\newcount\secno \global\secno=0
\global\newcount\meqno \global\meqno=1
\def\newsec#1{\global\advance\secno by1\message{(\the\secno. #1)}
\global\subsecno=0\eqnres@t\noindent{\bf\the\secno. #1}
\writetoca{{\secsym} {#1}}\par\nobreak\medskip\nobreak}
\def\eqnres@t{\xdef\secsym{\the\secno.}\global\meqno=1\bigbreak\bigskip}
\def\sequentialequations{\def\eqnres@t{\bigbreak}}\xdef\secsym{}
\global\newcount\subsecno \global\subsecno=0
\def\subsec#1{\global\advance\subsecno by1\message{(\secsym\the\subsecno. #1)}
\ifnum\lastpenalty>9000\else\bigbreak\fi
\noindent{\it\secsym\the\subsecno. #1}\writetoca{\string\quad 
{\secsym\the\subsecno.} {#1}}\par\nobreak\medskip\nobreak}
\def\appendix#1#2{\global\meqno=1\global\subsecno=0\xdef\secsym{\hbox{#1.}}
\bigbreak\bigskip\noindent{\bf Appendix #1. #2}\message{(#1. #2)}
\writetoca{Appendix {#1.} {#2}}\par\nobreak\medskip\nobreak}
%
%
\def\eqnn#1{\xdef #1{(\secsym\the\meqno)}\writedef{#1\leftbracket#1}%
\global\advance\meqno by1\wrlabeL#1}
\def\eqna#1{\xdef #1##1{\hbox{$(\secsym\the\meqno##1)$}}
\writedef{#1\numbersign1\leftbracket#1{\numbersign1}}%
\global\advance\meqno by1\wrlabeL{#1$\{\}$}}
\def\eqn#1#2{\xdef #1{(\secsym\the\meqno)}\writedef{#1\leftbracket#1}%
\global\advance\meqno by1$$#2\eqno#1\eqlabeL#1$$}
%
\newskip\footskip\footskip14pt plus 1pt minus 1pt 
\def\footnotefont{\ninepoint}\def\f@t#1{\footnotefont #1\@foot}
\def\f@@t{\baselineskip\footskip\bgroup\footnotefont\aftergroup\@foot\let\next}
\setbox\strutbox=\hbox{\vrule height9.5pt depth4.5pt width0pt}
\global\newcount\ftno \global\ftno=0
\def\foot{\global\advance\ftno by1\footnote{$^{\the\ftno}$}}
%
\newwrite\ftfile   
\def\footend{\def\foot{\global\advance\ftno by1\chardef\wfile=\ftfile
$^{\the\ftno}$\ifnum\ftno=1\immediate\openout\ftfile=foots.tmp\fi%
\immediate\write\ftfile{\noexpand\smallskip%
\noexpand\item{f\the\ftno:\ }\pctsign}\findarg}%
\def\footatend{\vfill\eject\immediate\closeout\ftfile{\parindent=20pt
\centerline{\bf Footnotes}\nobreak\bigskip\input foots.tmp }}}
\def\footatend{}
%
%
\global\newcount\refno \global\refno=1
\newwrite\rfile
\def\ref{[\the\refno]\nref}
\def\nref#1{\xdef#1{[\the\refno]}\writedef{#1\leftbracket#1}%
\ifnum\refno=1\immediate\openout\rfile=refs.tmp\fi
\global\advance\refno by1\chardef\wfile=\rfile\immediate
\write\rfile{\noexpand\item{#1\ }\reflabeL{#1\hskip.31in}\pctsign}\findarg}
\def\findarg#1#{\begingroup\obeylines\newlinechar=`\^^M\pass@rg}
{\obeylines\gdef\pass@rg#1{\writ@line\relax #1^^M\hbox{}^^M}%
\gdef\writ@line#1^^M{\expandafter\toks0\expandafter{\striprel@x #1}%
\edef\next{\the\toks0}\ifx\next\em@rk\let\next=\endgroup\else\ifx\next\empty%
\else\immediate\write\wfile{\the\toks0}\fi\let\next=\writ@line\fi\next\relax}}
\def\striprel@x#1{} \def\em@rk{\hbox{}} 
\def\lref{\begingroup\obeylines\lr@f}
\def\lr@f#1#2{\gdef#1{\ref#1{#2}}\endgroup\unskip}

\def\addref#1{\immediate\write\rfile{\noexpand\item{}#1}} 
\def\startrefs#1{\immediate\openout\rfile=refs.tmp\refno=#1}
\def\xref{\expandafter\xr@f}\def\xr@f[#1]{#1}
\def\refs#1{\count255=1[\r@fs #1{\hbox{}}]}
\def\r@fs#1{\ifx\und@fined#1\message{reflabel \string#1 is undefined.}%
\nref#1{need to supply reference \string#1.}\fi%
\vphantom{\hphantom{#1}}\edef\next{#1}\ifx\next\em@rk\def\next{}%
\else\ifx\next#1\ifodd\count255\relax\xref#1\count255=0\fi%
\else#1\count255=1\fi\let\next=\r@fs\fi\next}
%

%
\newwrite\ffile\global\newcount\figno \global\figno=1
\def\fig{fig.~\the\figno\nfig}
\def\nfig#1{\xdef#1{fig.~\the\figno}%
\writedef{#1\leftbracket fig.\noexpand~\the\figno}%
\ifnum\figno=1\immediate\openout\ffile=figs.tmp\fi\chardef\wfile=\ffile%
\immediate\write\ffile{\noexpand\medskip\noexpand\item{Fig.\ \the\figno. }
\reflabeL{#1\hskip.55in}\pctsign}\global\advance\figno by1\findarg}
\def\xfig{\expandafter\xf@g}\def\xf@g fig.\penalty\@M\ {}
\def\figs#1{figs.~\f@gs #1{\hbox{}}}
\def\f@gs#1{\edef\next{#1}\ifx\next\em@rk\def\next{}\else
\ifx\next#1\xfig #1\else#1\fi\let\next=\f@gs\fi\next}
\newwrite\lfile
{\escapechar-1\xdef\pctsign{\string\%}\xdef\leftbracket{\string\{}
\xdef\rightbracket{\string\}}\xdef\numbersign{\string\#}}

\def\writestop{\def\writestoppt{\immediate\write\lfile{\string\pageno%
\the\pageno\string\startrefs\leftbracket\the\refno\rightbracket%
\string\def\string\secsym\leftbracket\secsym\rightbracket%
\string\secno\the\secno\string\meqno\the\meqno}\immediate\closeout\lfile}}
\def\writestoppt{}\def\writedef#1{}
\def\seclab#1{\xdef #1{\the\secno}\writedef{#1\leftbracket#1}\wrlabeL{#1=#1}}
\def\subseclab#1{\xdef #1{\secsym\the\subsecno}%
\writedef{#1\leftbracket#1}\wrlabeL{#1=#1}}
\newwrite\tfile \def\writetoca#1{}
\def\leaderfill{\leaders\hbox to 1em{\hss.\hss}\hfill}
\def\writetoc{\immediate\openout\tfile=toc.tmp 
   \def\writetoca##1{{\edef\next{\write\tfile{\noindent ##1 
   \string\leaderfill {\noexpand\number\pageno} \par}}\next}}}

%
\catcode`\@=12 
%
\edef\tfontsize{\ifx\answ\bigans scaled\magstep3\else scaled\magstep4\fi}
\font\titlerm=cmr10 \tfontsize \font\titlerms=cmr7 \tfontsize
\font\titlermss=cmr5 \tfontsize \font\titlei=cmmi10 \tfontsize
\font\titleis=cmmi7 \tfontsize \font\titleiss=cmmi5 \tfontsize
\font\titlesy=cmsy10 \tfontsize \font\titlesys=cmsy7 \tfontsize
\font\titlesyss=cmsy5 \tfontsize \font\titleit=cmti10 \tfontsize
\skewchar\titlei='177 \skewchar\titleis='177 \skewchar\titleiss='177
\skewchar\titlesy='60 \skewchar\titlesys='60 \skewchar\titlesyss='60
\def\titlefont{\def\rm{\fam0\titlerm}
\textfont0=\titlerm \scriptfont0=\titlerms \scriptscriptfont0=\titlermss
\textfont1=\titlei \scriptfont1=\titleis \scriptscriptfont1=\titleiss
\textfont2=\titlesy \scriptfont2=\titlesys \scriptscriptfont2=\titlesyss
\textfont\itfam=\titleit \def\it{\fam\itfam\titleit}\rm}
 \ifx\answ\bigans\else scaled\magstep1\fi
\ifx\answ\bigans\def\abstractfont{\tenpoint}\else
\font\abssl=cmsl10 scaled \magstep1
\font\absrm=cmr10 scaled\magstep1 \font\absrms=cmr7 scaled\magstep1
\font\absrmss=cmr5 scaled\magstep1 \font\absi=cmmi10 scaled\magstep1
\font\absis=cmmi7 scaled\magstep1 \font\absiss=cmmi5 scaled\magstep1
\font\abssy=cmsy10 scaled\magstep1 \font\abssys=cmsy7 scaled\magstep1
\font\abssyss=cmsy5 scaled\magstep1 \font\absbf=cmbx10 scaled\magstep1
\skewchar\absi='177 \skewchar\absis='177 \skewchar\absiss='177
\skewchar\abssy='60 \skewchar\abssys='60 \skewchar\abssyss='60
\def\abstractfont{\def\rm{\fam0\absrm}
\textfont0=\absrm \scriptfont0=\absrms \scriptscriptfont0=\absrmss
\textfont1=\absi \scriptfont1=\absis \scriptscriptfont1=\absiss
\textfont2=\abssy \scriptfont2=\abssys \scriptscriptfont2=\abssyss
\textfont\itfam=\bigit \def\it{\fam\itfam\bigit}\def\footnotefont{\tenpoint}%
\textfont\slfam=\abssl \def\sl{\fam\slfam\abssl}%
\textfont\bffam=\absbf \def\bf{\fam\bffam\absbf}\rm}\fi
\def\tenpoint{\def\rm{\fam0\tenrm}
\textfont0=\tenrm \scriptfont0=\sevenrm \scriptscriptfont0=\fiverm
\textfont1=\teni  \scriptfont1=\seveni  \scriptscriptfont1=\fivei
\textfont2=\tensy \scriptfont2=\sevensy \scriptscriptfont2=\fivesy
\textfont\itfam=\tenit \def\it{\fam\itfam\tenit}\def\footnotefont{\ninepoint}%
\textfont\bffam=\tenbf \def\bf{\fam\bffam\tenbf}\def\sl{\fam\slfam\tensl}\rm}
\font\ninerm=cmr9 \font\sixrm=cmr6 \font\ninei=cmmi9 \font\sixi=cmmi6 
\font\ninesy=cmsy9 \font\sixsy=cmsy6 \font\ninebf=cmbx9 
\font\nineit=cmti9 \font\ninesl=cmsl9 \skewchar\ninei='177
\skewchar\sixi='177 \skewchar\ninesy='60 \skewchar\sixsy='60 
\def\ninepoint{\def\rm{\fam0\ninerm}
\textfont0=\ninerm \scriptfont0=\sixrm \scriptscriptfont0=\fiverm
\textfont1=\ninei \scriptfont1=\sixi \scriptscriptfont1=\fivei
\textfont2=\ninesy \scriptfont2=\sixsy \scriptscriptfont2=\fivesy
\textfont\itfam=\ninei \def\it{\fam\itfam\nineit}\def\sl{\fam\slfam\ninesl}%
\textfont\bffam=\ninebf \def\bf{\fam\bffam\ninebf}\rm} 
%
%
\def\noblackbox{\overfullrule=0pt}
\hyphenation{anom-aly anom-alies coun-ter-term coun-ter-terms}
\def\inv{^{\raise.15ex\hbox{${\scriptscriptstyle -}$}\kern-.05em 1}}

\def\Dsl{\,\raise.15ex\hbox{/}\mkern-13.5mu D} 
\def\dsl{\raise.15ex\hbox{/}\kern-.57em\partial}

\font\bigit=cmti10 scaled \magstep1
\def\lspace{\ifx\answ\bigans{}\else\qquad\fi}
\def\lbspace{\ifx\answ\bigans{}\else\hskip-.2in\fi} 
\def\boxeqn#1{\vcenter{\vbox{\hrule\hbox{\vrule\kern3pt\vbox{\kern3pt
	\hbox{${\displaystyle #1}$}\kern3pt}\kern3pt\vrule}\hrule}}}
\def\mbox#1#2{\vcenter{\hrule \hbox{\vrule height#2in
		\kern#1in \vrule} \hrule}}  
%
    
  \def\CI{{\cal I}}

\def\darr#1{\raise1.5ex\hbox{$\leftrightarrow$}\mkern-16.5mu #1}

\def\roughly#1{\raise.3ex\hbox{$#1$\kern-.75em\lower1ex\hbox{$\sim$}}}

\def\frac#1#2{{#1\over#2}}

\def\journal#1&#2(#3){\unskip, #1~\bf #2 \rm(19#3) }
\def\andjournal#1&#2(#3){\sl #1~\bf #2 \rm (19#3) }

\def\bra#1{\left\langle #1\right|}
\def\ket#1{\left| #1\right\rangle}

\catcode`\@=11\def\slash#1{\mathord{\mathpalette\c@ncel{#1}}}
\overfullrule=0pt
\def\steepslash{\c@ncel}
\def\frac#1#2{{#1\over #2}}

\def\:{\!:\!}
\def\inbar{\,\vrule height1.5ex width.4pt depth0pt}
\def\IQ{\relax\,\hbox{$\inbar\kern-.3em{\rm Q}$}}
\def\IB{\relax{\rm I\kern-.18em B}}
\def\IC{\relax\hbox{$\inbar\kern-.3em{\rm C}$}}
\def\IP{\relax{\rm I\kern-.18em P}}
\def\IR{\relax{\rm I\kern-.18em R}}
\def\ZZ{\relax\ifmmode\mathchoice
{\hbox{Z\kern-.4em Z}}{\hbox{Z\kern-.4em Z}}
{\lower.9pt\hbox{Z\kern-.4em Z}}
{\lower1.2pt\hbox{Z\kern-.4em Z}}\else{Z\kern-.4em Z}\fi}

\catcode`\@=12

\def\npb#1(#2)#3{{ Nucl. Phys. }{B#1} (#2) #3}
\def\plb#1(#2)#3{{ Phys. Lett. }{#1B} (#2) #3}
\def\pla#1(#2)#3{{ Phys. Lett. }{#1A} (#2) #3}
\def\prl#1(#2)#3{{ Phys. Rev. Lett. }{#1} (#2) #3}
\def\mpla#1(#2)#3{{ Mod. Phys. Lett. }{A#1} (#2) #3}
\def\ijmpa#1(#2)#3{{ Int. J. Mod. Phys. }{A#1} (#2) #3}
\def\cmp#1(#2)#3{{ Comm. Math. Phys. }{#1} (#2) #3}
\def\cqg#1(#2)#3{{ Class. Quantum Grav. }{#1} (#2) #3}
\def\jmp#1(#2)#3{{ J. Math. Phys. }{#1} (#2) #3}
\def\anp#1(#2)#3{{ Ann. Phys. }{#1} (#2) #3}
\def\prd#1(#2)#3{{ Phys. Rev. } {D{#1}} (#2) #3}
\def\ptp#1(#2)#3{{ Progr. Theor. Phys. }{#1} (#2) #3}
\def\aom#1(#2)#3{{ Ann. Math. }{#1} (#2) #3}

\def\bra{\langle}
\def\ket{\rangle}

\def\C{{\bf C}}
\def\R{{\bf R}}
\def\Z{{\bf Z}}
\def\P{{\bf P}}

\def\Z{{\bf Z}}
\def\cA{{\cal A}}
\def\cG{{\cal G}}
\def\cI{{\cal I}}
\def\cJ{{\cal J}}

\def\cM{{\cal M}}

\def\cC{{\cal C}}
\def\cI{{\cal I}}

\def\cL{{\cal L}}

\def\cP{{\cal P}}
\def\cX{{\cal X}}

\def\cicy#1(#2|#3)#4{\left(\matrix{#2}\right|\!\!
                     \left|\matrix{#3}\right)^{{#4}}_{#1}}

\def\ra{\rightarrow}

\def\Box{{\,\lower0.9pt\vbox{\hrule 
\hbox{\vrule height 0.2 cm \hskip 0.2 cm  
\vrule height 0.2 cm}\hrule}\,}}

\global\newcount\thmno \global\thmno=0
\def\definition#1{\global\advance\thmno by1
\bigskip\noindent{\bf Definition \secsym\the\thmno. }{\it #1}
\par\nobreak\medskip\nobreak}
\def\question#1{\global\advance\thmno by1
\bigskip\noindent{\bf Question \secsym\the\thmno. }{\it #1}
\par\nobreak\medskip\nobreak}
\def\theorem#1{\global\advance\thmno by1
\bigskip\noindent{\bf Theorem \secsym\the\thmno. }{\it #1}
\par\nobreak\medskip\nobreak}
\def\proposition#1{\global\advance\thmno by1
\bigskip\noindent{\bf Proposition \secsym\the\thmno. }{\it #1}
\par\nobreak\medskip\nobreak}
\def\corollary#1{\global\advance\thmno by1
\bigskip\noindent{\bf Corollary \secsym\the\thmno. }{\it #1}
\par\nobreak\medskip\nobreak}
\def\lemma#1{\global\advance\thmno by1
\bigskip\noindent{\bf Lemma \secsym\the\thmno. }{\it #1}
\par\nobreak\medskip\nobreak}
\def\conjecture#1{\global\advance\thmno by1
\bigskip\noindent{\bf Conjecture \secsym\the\thmno. }{\it #1}
\par\nobreak\medskip\nobreak}
\def\exercise#1{\global\advance\thmno by1
\bigskip\noindent{\bf Exercise \secsym\the\thmno. }{\it #1}
\par\nobreak\medskip\nobreak}
\def\remark#1{\global\advance\thmno by1
\bigskip\noindent{\bf Remark \secsym\the\thmno. }{\it #1}
\par\nobreak\medskip\nobreak}
\def\problem#1{\global\advance\thmno by1
\bigskip\noindent{\bf Problem \secsym\the\thmno. }{\it #1}
\par\nobreak\medskip\nobreak}
\def\others#1#2{\global\advance\thmno by1
\bigskip\noindent{\bf #1 \secsym\the\thmno. }{\it #2}
\par\nobreak\medskip\nobreak}
\def\proof{\noindent Proof: }

\def\thmlab#1{\xdef #1{\secsym\the\thmno}\writedef{#1\leftbracket#1}
              \wrlabeL{#1=#1}}
\def\remlab#1{\xdef #1{\secsym\the\thmno}\writedef{#1\leftbracket#1}
              \wrlabeL{#1=#1}}

%
%
\def\newsec#1{\global\advance\secno by1\message{(\the\secno. #1)}
\global\subsecno=0\thmno=0\eqnres@t\noindent{\bf\the\secno. #1}
\writetoca{{\secsym} {#1}}\par\nobreak\medskip\nobreak}
\def\eqnres@t{\xdef\secsym{\the\secno.}\global\meqno=1\bigbreak\bigskip}
\def\sequentialequations{\def\eqnres@t{\bigbreak}}\xdef\secsym{}
\catcode`@=12

\baselineskip16pt
\noblackbox

\def\A{{\cal A}}

\def\IP{{\bf P}}
\def\IR{{\bf R}}
\def\IC{{\bf C}}
\def\IQ{{\bf Q}}
\def\ZZ{{\bf Z}}
\def\pd{\partial}

\def\({ \left(  }
\def\){ \right) }

\def\ta#1{\theta_{a_#1}}

\def\X#1#2(#3)#4#5{ {$X_{#1}^#2(#3)_{#5}^{#4}$} }

\def\batyrevJA{[1]}
 \def\refbatyrevJA{V. Batyrev,J. Algebraic Geometry 3 (1994) 493.}

\def\batyrevDuke{[2]}
 \def\refbatyrevDuke{V. Batyrev, Duke Math.J.{\bf 69} (1993) 349.}

\def\batyrevQ{[3]}
 \def\refbatyrevQ{V. Batyrev,{\sl Quantum cohomology rings of toric manifolds},
                   preprint 1993.}
\def\BKK{[4]}
 \def\refBKK{P.Berglund, S.Katz and A.Klemm, Nucl.Phys.{\bf B456}(1995)153.}

\def\BFS{[5]}
 \def\refBFS{L.Billera,P.Filliman and B.Sturmfels,
             Adv.in Math.{\bf 83}(1990)155.}
\def\candelasETAL{[6]}
 \def\refcandelasETAL{P.Candelas, X. de la Ossa, P.Green and L.Parks, 
                      Nucl.Phys.{\bf B359}(1991)21.}
\def\candelasIIa{[7]}
 \def\refcandelasIIa{P.Candelas, X. de la Ossa, A.Font, S.Katz and D.Morrison, 
                     Nucl.Phys.{\bf B416}(1994)481.}
\def\candelasIIb{[8]}
 \def\refcandelasIIb{P.Candelas, X. de la Ossa, A.Font, S.Katz and D.Morrison, 
                     Nucl.Phys.{\bf B429}(1994)629.}

\def\CoxLittleOshea{[9]}
 \def\refCoxLittleOshea{D. Cox, J. Little and D. O'Shea, 
            {\sl Ideals, Varieties and Algorithms}, UTM Springer-Verlag 1992.}
\def\Ft{[10]}
  \def\refFt{A.Font, Nucl.Phys.{\bf B391}(1993)358.}

\def\fulton{[11]}
 \def\reffulton{W. Fulton, {\sl An introduction to toric varieties}, 
                            Princeton Univ. Press 1993.}
\def\GKZ{[12]}
 \def\refGKZ{ I.M.Gel'fand, A.V.Zelevinsky and M.M.Kapranov, 
              Funct. Anal. Appl.{\bf 28}(1989)94. }
\def\GKZII{[13]}
 \def\refGKZII{ I.M.Gel'fand, A.V.Zelevinsky and M.M.Kapranov, 
              Adv. Math. {\bf 84}(1990). }
\def\GreenePlesser{[14]}
  \def\refGreenePlesser{B.Greene and M.Plesser, Nucl.Phys.{\bf B338}(1990)15.}

\def\HKTYI{[15]}
 \def\refHKTYI{S.Hosono, A.Klemm, S.Theisen and S.-T.Yau, 
               Commun. Math. Phys. 167 (1995) 301.}
\def\HKTYII{[16]}
 \def\refHKTYII{S.Hosono, A.Klemm, S.Theisen and S.-T.Yau, 
                Nucl. Phys. B433 (1995) 501.}
\def\HLY{[17]}
 \def\refHLY{S. Hosono, B. Lian and  S.T. Yau, 
      {\sl GKZ-Generalized Hypergeometric Systems 
           in Mirror Symmetry of Calabi-Yau Hypersurfaces}, 
            Harvard Univ. preprint, alg-geom/9511001, to appear in CMP 1996.}
\def\KT{[18]}
 \def\refKT{A.Klemm and S.Theisen, Nucl.Phys.{\bf B389}(1993)153.}

\def\morrison{[19]}
 \def\refmorrison{D.Morrison, 
      {\it Picard-Fuchs Equations and Mirror Maps For Hypersurfaces},
      in {\it Essays on Mirror Manifolds}, 
      Ed. S.-T.Yau, International Press 1992.}

\def\MFK{[20]}
 \def\refMFK{D. Mumford, J. Fogarty and F. Kirwan, 
             {\sl Geometric Invariant Theory}, 3rd Ed., Springer-Verlag 1994.} 
\def\oda{[21]}
 \def\refoda{T. Oda, {\sl Convex Bodies and Algebraic Geometry}, 
                     Springer-Verlag 1988.} 
\def\OdaPark{[22]}
 \def\refOdaPark{T. Oda and H.S. Park, T\^ohoku Math. J. {\bf 43}(1991) 375-399.}

\def\sturmfels{[23]}
 \def\refsturmfels{B.Sturmfels,T\^ohoku Math. J. {\bf 43} (1991) 249.} 

\def\sturmfelsII{[24]}
 \def\refsturmfelsII{B.Sturmfels, {\sl Gr\"obner Bases and Convex
Polytopes}, AMS University Lecture Series, Vol. 8, Providence RI, 1995. } 

\def\yau{[25]}
 \def\refyau{{\it Essays on Mirror Manifolds}, 
      Ed. S.-T.Yau, International Press 1992.}


\pageno=1

\Title{alg-geom/9603014}
{\vbox{
 \centerline{Maximal Degeneracy Points of GKZ Systems}
    } }

   \centerline{S. Hosono$^1$\footnote{$^\dagger$}
   {email: hosono@sci.toyama-u.ac.jp},
   B.H. Lian$^2$\footnote{$^\ddagger$}{email: lian@max.math.brandeis.edu}
   and S.-T. Yau$^3$\footnote{$^\diamond$}{email: yau@math.harvard.edu}}

   \bigskip\centerline{
   \vbox{
   \hbox{ $^1$ Department of Mathematics}
   \hbox{ \hskip30pt Toyama University}
   \hbox{ \hskip28pt Toyama 930, Japan} }
   \hskip0.5cm
   \vbox{
   \hbox{ $^2$ Department of Mathematics}
   \hbox{ \hskip30pt Brandeis University}
   \hbox{ \hskip27pt Waltham, MA 02154} } }
   \bigskip
   \centerline{
   \vbox{
   \hbox{ $^3$ Department of Mathematics}
   \hbox{ \hskip30pt Harvard University}
   \hbox{ \hskip23pt Cambridge, MA 02138} } }

\vskip.2in
Abstract. Motivated by mirror symmetry, we study certain integral
representations of solutions to the Gel'fand-Kapranov-Zelevinsky(GKZ) 
hypergeometric system.
Some of these solutions arise as period integrals for Calabi-Yau manifolds
in mirror symmetry.
We prove that for a suitable compactification
of the parameter space, there exists certain special boundary points,
which we called maximal degeneracy points, at which all but
one solutions become singular. 

\Date{3/2/96} 

\newsec{Introduction}

In this paper, we study the  
Gel'fand-Kapranov-Zelevinsky hypergeometric PDE systems{\GKZ}.
We are mainly concerned with a case of the so-called resonant
exponent $\beta$, though some of our results apply to generic $\beta$
as well. This resonant case in fact plays an especially important
role for applications in mirror symmetry.

The problem we consider is motivated by mirror symmetry{\yau} as follows.
Given a family of Calabi-Yau manifolds $\pi:\cX\ra S$, one would
like to give a suitable compactification in which there is a boundary  
point, called a large radius limit, where all but one of period integrals
of the manifolds become singular near that point. 
An astonishing discovery{\candelasETAL} in mirror symmetry is that
under some assumptions, the period integrals near this point
combine to form a generating function which predicts the
number of rational curves (a.k.a. instanton contributions)
on another Calabi-Yau manifold.
Thus the existence and construction of 
 a large radius limit
pose an important problem.

When the family of Calabi-Yau manifolds are anticanonical divisors
in a suitable toric variety, then their period integrals are
natural solutions to a GKZ system{\batyrevDuke}. 
Thus the GKZ theory becomes an
important tool for studying the behavior of such period integrals
{\HKTYI}{\HKTYII}{\BKK}{\HLY}.
Motivated by this, we study a GKZ system canonically
associated with a complete regular fan, and
 a class of integral solutions. We prove a general existence
theorem for the so-called maximal degeneracy points{\morrison} for these
integral representations. These special points can be geometrically
interpreted as fixed points under a canonical torus action on
the parameter space. Some complete regular fans (in most interesting
examples we know) arise in mirror symmetry{\batyrevJA}. When
they do, then the maximal degeneracy points turn out to be
large radius limits where the enumeration of instanton
contributions is carried out\candelasETAL\morrison
{\KT}{\Ft}{\candelasIIa}{\HKTYI}{\candelasIIb}{\HKTYII}{\BKK}{\HLY}.
 
\vfill\eject

\subsec{Conventions and Notations}

\item{$N$}: rank $n$ lattice
\item{$M$}: dual lattice
\item{$N_\R$}: $N\otimes\R$
\item{$M_\R$}: $M\otimes\R$
\item{$\bar{N}$}: $\Z\times N$
\item{$\bar{M}$}: $\Z\times M$
\item{$\mu$}: a point in $N$
\item{$\bar\mu$}: $1\times\mu$
\item{$\Sigma$}: a complete rational polyhedral fan in $N$
\item{$\cA$}: a finite set of integral points which spans
              the affine hyperplane $1\times N_\R\subset\bar N_\R$
\item{$L=L_\cA$}: the lattice of relations on $\cA$
\item{$\P_\Sigma$}: the toric variety associated to the fan $\Sigma$
\item{$\Sigma(1)$}: the set of primitive elements on the edges of the fan 
$\Sigma$
\item{$\sigma(1)$}: the set of primitive elements on the edges of the cone
$\sigma$
\item{$S\Sigma$}: the secondary fan of $\Sigma$
\item{$G\Sigma$}: the Gr\"obner fan
\item{$cone(S)$}: the cone generated by a finite set $S$
\item{$conv(S)$}: the convex hull of a finite set $S$

\bigskip
An $n$ dimensional cone in $N_\R$ is also called a large cone. A cone
which contains no nontrivial linear subspace is called strongly
convex. A cone generated by finitely many integral points are called
rational. A rational cone is called regular if it is generated
by a subset of an integral base. A fan is regular if it consists of 
regular cones. The dual of a cone $\sigma$ is denoted $\sigma^\vee$.

\vfill\eject

\newsec{GKZ $\cA$-hypergeometric systems}

Let $\cA\subset\bar N$ be a set of $p+1$ integral points which span the
affine plane $1\times N_\R$. We sometimes write
$\cA=\{\bar\mu_0,..,\bar\mu_p\}$. The projection map $\bar N\ra N$ is
denoted $\bar\mu\mapsto\mu$.
Let $L=L_\cA:=\{l\in\Z^{p+1}|\sum_{\bar\mu\in\cA} l_\mu\bar\mu=0\}$, 
which we call the lattice of relations on $\cA$. 
This is clearly a lattice of rank $p-n$.
The GKZ system with exponent
$\beta\in\bar{N}_\R$ is the PDE system on $\C^{p+1}=\C^\cA$ given by
\eqn\dumb{\Box_l\Pi(a)=0,~~~\left(\sum\bra
u,\bar{\mu}\ket\theta_{a_\mu}-\bra
u,\beta\ket\right)\Pi(a)=0}
for $l\in L$, $u\in \bar{M}_\R${\GKZ}.
Here $\Box_l:=\prod_{l_\mu>0}({\partial\over\partial a_\mu})^{l_\mu}
-\prod_{l_\mu<0}({\partial\over\partial a_\mu})^{-l_\mu}$
and $\theta_{a_\mu}:=a_\mu{\partial~\over\partial a_\mu}$.
Equivalently we can write
\eqn\dumb{a^{-\gamma}\Box_l a^\gamma f(a)=0,~~~\sum\bra
u,\bar{\mu}\ket\theta_{a
_\mu}f(a)=0}
where $f(a)=a^{-\gamma}\Pi(a)$ and
$\gamma$ lies in the $(p-n)$-dimensional affine subspace
$\Phi(\beta):=\{\gamma
|\sum_\mu\gamma_\mu\bar\mu=\beta\}$.
{\it We shall always assume that 
$\bar\mu_0:=1\times\vec0\in\cA$.}

\newsec{Period Integrals}

We shall consider a class of (multi-valued) analytic functions
attached to the data $\cA$ as follows. Given
a choice of basis $e_1,..,e_n$ 
of $N$, we have an isomorphism $T_M=Hom(N,\C^\times)\cong
(\C^\times)^n$. On this algebraic torus, there is a family of rational
functions defined by the 
Laurent polynomials:
\eqn\dumb{f_\cA(X,a)=\sum_{\bar\mu\in\cA} a_\mu X^\mu.}
The zero locus of this function is a divisor $H_a$
in the torus $T_M$. Now $T_M$ has a canonical $T_M$-invariant
holomorphic $n$-form given by $\prod {dX_i\over X_i}$. Following
{\batyrevDuke} for
every homology $n$-cycle $\alpha$ in the manifold $T_M-H_a$,
we define the {\it period integral}
\eqn\dumb{\Pi_\alpha(a)=\int_\alpha{1\over f_\cA(X,a)}\prod
{dX_i\over X_i}.}
This integral converges in some domain, and is singular in
some divisor $D\subset\C^\cA$. Under analytic
continuation in $\C^\cA-D$,
the integral will in general have monodromy. Thus the period integrals are
local sections of a locally constant sheaf over $\C^\cA-D$.
As an example, we
have the following. Let $\alpha_0$ be the real torus cycle
$\{|X_1|=\cdots=|X_n|=1\}$. 
The following two results are proved in {\batyrevDuke} 

\proposition{
In the domain $|a_0|>>|a_1|,..,|a_p|$, we have the convergent
power
 series representation
 \eqn\dumb{a_0\int_{\alpha_0}{1\over f_\cA(X,a)}
 \prod_{i=1}^n{dX_i\over X_i}=
 (2\pi i)^n \sum_{l_1\mu_1+\cdots+l_p\mu_p=0,~l_1,..,l_p\geq0}
(-1)^{l_1+\cdots+l_p}
{(l_1+\cdots+l_p)!\over l_1!\cdots l_p!}
 {a_1^{l_1}\cdots a_p^{l_p}\over a_0^{l_1+\cdots+l_p}}.} }
 \thmlab\PowerSeries
 \proof First we have the geometric series ${a_0\over f_\cA(X,a)}=
 \sum_{r\geq0}{1\over a_0^r}(-a_1X^{\mu_1}-\cdots-a_pX^{\mu_p})^r$. Also we
 have
 $\int_{\alpha_0}X^\lambda\prod_{i=1}^n{dX_i\over X_i}=
(2\pi i)^n \delta_{\lambda,0}$. Thus the terms in the geometric series
 which contribute to the integral are powers of the form
$(-1)^{l_1+\cdots+l_p}
{a_1^{l_1}\cdots a_p^{l_p}\over a_0^{l_1+\cdots+l_p}}$ where
 $l_1\mu_1+\cdots+l_p\mu_p=0,~l_1,..,l_p\geq0$. In each
 term  ${1\over a_0^r}(-a_1X^{\mu_1}-\cdots-a_pX^{\mu_p})^r$, each
 tuple $(l_1,..,l_p)$ with $r=l_1+\cdots+l_p$ contributes
 exactly ${(l_1+\cdots+l_p)!\over l_1!\cdots l_p!}$ copies of such a power.
 Now summing over all $r\geq0$ and over all those tuples, we get
 the desired sum. The sum is convergent because it is a subseries of
 the standard series (with $z_i=-a_i/a_0$):
 \eqn\dumb{\sum_{l_1,..,l_p\geq0}
 {(l_1+\cdots+l_p)!\over l_1!\cdots l_p!}z_1^{l_1}\cdots z_p^{l_p}. ~~~\Box}

\proposition{(see {\GKZII}) The period integrals $\Pi_\alpha(a)$
are solutions to the GKZ system
with the exponent $\beta=-1\times\vec0$.}
\proof We claim that $\Box_l\Pi_\alpha(a)=0$ for all $l\in L$. In fact,
$\Box_lG(f_\cA(a,X))=0$
for any differentiable function $G$. First we have
\eqn\dumb{\partial_{a_{j_1}}\cdots\partial_{a_{j_k}}G(f_\cA(a,X))=G^{(k)}(f)
X^{\mu_{j_1}}\cdots X^{\mu_{j_k}} }
Write $l\in L$ as $l=l^+-l^-$ where $l^\pm$ have only nonnegative entries
and have disjoint support. Note that the sum of the entries of each $l^+$
is the same as that of $l^-$. Call this sum $k$. Then we have
\eqn\dumb{\Box_lG(f_\cA(a,X))=G^{(k)}(f_\cA)
(X^{\sum l_{\mu}^+ \mu}-X^{\sum l_{\mu}^- \mu}). }
But $l\in L$ implies that $0=\sum l_\mu\mu=\sum (l_{\mu}^+-l_{\mu}^-)\mu$
It follows that $\Box_lG(f_\cA(a,X))=0$.

Now we have to check that the periods $\Pi_\alpha(a)$
are annihilated by the order
1 operators in the GKZ system.
We consider the action of $\lambda:=(\lambda_0,\lambda_1,..,\lambda_n)$
on the periods as follows: under $X_i\mapsto\lambda_i X_i$, the cycle $\alpha$
gets moved to $\alpha_\lambda$ and $dm:=\prod^n_{i=1}{dX_i\over X_i}$ is
invariant. Thus we get
\eqn\dumb{{1\over\lambda_0}\int_{\alpha}{1\over \sum a_\mu X^\mu} dm
=\int_{\alpha_\lambda}{1\over \sum a_\mu\lambda^{\bar\mu}X^\mu} dm.}
Now for $\lambda$ closed enough to the identity element, it's clear that
the cycle $\alpha_\lambda$ is homologous to $\alpha$. So we can replace
$\alpha_\lambda$ on the RHS by $\alpha$. Differentiating the equation above
by $\lambda_i{\partial\over\partial\lambda_i}$ 
and then set $\lambda=id=(1,..,1)$, we get
\eqn\dumb{\eqalign{
-{1\over\lambda_0}|_{\lambda=id}
\delta_{i,0}\int_{\alpha}{1\over \sum a_\mu X^\mu} dm
&=\int_{\alpha}{-\sum a_\mu\bar\mu_i\lambda^{\bar\mu}X^\mu\over (\sum
a_\mu\lambda^{\bar\mu}X^\mu)^2} dm|_{\lambda=id}\cr
&=\int_{\alpha}{-\sum \bra e_i,\bar\mu\ket a_\mu X^\mu\over
f_\cA(a,X)^2} dm\cr
&=\sum \bra e_i,\bar\mu\ket\theta_{a_\mu}\int_{\alpha}{1\over f_\cA(a,X)}dm} }
where $\theta_{a_\mu}:=a_\mu{\partial\over\partial a_\mu}$.
It follows that $\bra e_i,\beta\ket\Pi_\alpha(a)
=\sum_\mu \bra e_i,\bar\mu\ket\theta_{a_\mu}\Pi_\alpha(a)$. $\Box$

\remark{Note that the order 1 operators in the GKZ system are
given by the vector field of the algebraic group $T_M$ with some twist
given by the exponent $\beta=-1\times\vec0$. Thus the order 1 equations 
means that
the period integrals are of the form $a_0^{-1}$ times a function
invariant under $\C^\times\times T_M$ 
which acts by $\lambda:(a_\mu)\mapsto(a_\mu\lambda^{\bar\mu})$, where
$\lambda=(\lambda_0,..,\lambda_n)$.}
\thmlab\remactions

\subsec{the maximal degeneracy problem}

Since the GKZ system is defined on a domain
in $\C^{p+1}=\C^\cA$, it is natural to study how the solutions degenerate
near the ``boundary'' of the domain. This should tell us something
about how the period integrals of the hypersurfaces $f_\cA=0$
degenerates. This of course depends on how we compactify
the domain. Now any nonzero $f_\cA$ defines a hypersurface.
Thus it's natural to consider the projective space $\P\C^{\cA}$.
As remarked in {\remactions} the GKZ system has a symmetry group given by the
algebraic torus $\C^\times\times T_M\cong(\C^\times)^{n+1}$.
It seems natural to consider the quotient 
$\P\C^{\cA}/(\C^\times\times T_M)$. 
If $G$ is an algebraic group acting on a variety $V$, $V/G$ in general
won't be a variety. However
note that $\P\C^{\cA}$ admits a canonical action by $(\C^\times)^{\cA}$,
and the group $\C^\times\times T_M$ acts by a subgroup.
Thus we expect ``$\P\C^\cA/(\C^\times\times T_M)$" to be a toric variety with
a $(p-n)$-dimensional open torus.
There are two ways to repair the quotient. We can  use geometric invariant
theory {\MFK} to construct a model for the quotient. We shall however use a
different approach which we now explain. A virtue of this approach is
that the resulting compactification can be described quite explicitly.

\proposition{The map $(\C^\times)^\cA/(\C^\times\times T_M)\ra
Hom(L,\C^\times)$ given by
$(a_\mu)\mapsto x_a$
where $x_a(l)=(-1)^{l_0}a^l$ is an isomorphism.}
\thmlab\PrincipalBundle

The extra sign is only a matter of convenience here, but is important
in applications to mirror symmetry in physics{\HKTYI}.
We denote by $S_\cA$ the image of $(\C^\times)^\cA-D$ under the
above isomorphism. This is the complement of a divisor in $Hom(L,\C^\times)$.
By the standard construction of toric varieties from fans, we see
that
any toric variety $\P_F$ associate with a complete fan $F$ in $L^*_\R$
is a compactification of
$Hom(L,\C^\times)$ and hence of $S_\cA$.
Since the period integrals
$\bar\Pi_\alpha:=a_0\Pi_\alpha(a)$ are $\C^\times\times T_M$-invariant,
they descend to local sections of a locally constant sheaf over $S_\cA$.

\definition{We call a smooth boundary point $p\in\P_F-Hom(L,\C^\times)$
a maximal degeneracy point if near $p$ there is exactly one period integral (up
to multiple) $\bar\Pi_\alpha$ which extends holomorphically
across $p$.}

 We will find in section 6 a natural
compactification $\P_F$ -- 
one which can be described in terms of the initial data $\cA$
alone. And we will prove that the existence of 
maximal degeneracy points in this compactification.

\newsec{The Secondary Fan and Gr\"obner Fan}

\subsec{the secondary fan}

Let us consider the convex polytope $conv(\cA)$. 
A triangulation of $conv(\cA)$ (see {\GKZ}) is a collection of
$n$-simplices whose vertices are in $\cA$, such that
the intersection of two such simplices is a face of both.
Thus $T$ is given by a collection of bases $I\subset\cA$ of $\bar N_\R$,
each basis being the set of vertices of an $n$-simplex.
A continuous function $h$
on the cone over $conv(\cA)\subset \bar{N}_\R$ is $T$-piecewise
linear if it is linear in the cone over each $n$-simplex of $T$. 
Note that $h$ is determined by its values at the vertices $\bar{\mu}$
of $n$-simplices in $T$. Thus each point
in $x\in\R^{p+1}=\R^\cA$ 
determines a piecewise linear function $h_x$ whose values
at a vertex $\bar\mu_i$ is $x_i$. 
A piecewise linear function $h_x$ is called convex 
if $h_x(a+b)\leq h_x(a)+h_x(b)$ for arbitrary a $a,b$ in 
the cones, and is called strictly convex if in addition 
it satisfies $h_x \vert_\sigma \not= h_x \vert_\tau$  
for any large cones $\sigma, \tau$. 

Let $\cC'(T)$
be the set of points $x$ for which $h_x$ is convex 
and has the property that $h_x(\bar\mu_i)\leq x_i$ if 
$\bar\mu_i$ is not a vertex of $T$. 
It is known that  
$\cC'(T)$ is a polyhedral cone in $\R^{p+1}$ but 
it is not strongly convex. 
The triangulation $T$ is called regular
if $\cC'(T)$ has an interior point, i.e., 
there is a strictly convex $T$-piecewise
linear function.
Define for each basis $I$ of $N_\R$, 
\eqn\CI{\cC'(I):=\{x\in
\R^{p+1}|h_{I,x}(\bar\mu_j)\leq x_j,~\forall j\notin I\}}
where $h_{I,x}$ is the linear function with value $x_i$ at $\bar\mu_i$,
$i\in I$. We also define for the basis $I$ a subspace of $L_\R$, 
\eqn\KI{ K(I):=\{x\in L_R|x_i\geq0,~\forall i\notin I\}. } 
Then  

\proposition{{\GKZ} 
If $I\subset\cA$ is a basis of $\bar N_\R$, then $\cC'(I)^\vee=K(I)$. 
If $T$ is a triangulation of $conv(\cA)$, then
$\cC'(T)=\cap_{I\in T} \cC'(I)$.} 
\thmlab\GKZprop

The collection of the cones $\cC'(T)$ with $T$ regular, 
together with their faces form a fan
called the secondary fan $S\Sigma$. 
Now every cone contains a minimal face which
is the a linear space of dimension $n+1$. This is exactly the space 
of points $x$ for which $\{h_x\vert_\sigma\}$ glues to a global 
linear function. 
This space may be naturally identified with $\bar{N}_\R^*=\bar{M}_\R$.
If we project the secondary fan 
in $\R^{p+1}$ along
the $\bar{M}_\R$, we obtain a new fan in $\R^{p+1}/\bar{M}_\R$. We also call
this the secondary fan. It is known that the secondary fan now consists of 
strongly convex cones and is a complete fan. 
The image of the
cone $\cC'(T)$ under the projection will be denoted $\cC(T)$.
Again under the natural identification
$\R^{p+1}/\bar{M}_\R=L_\R^*$,
we can regard each cone $\cC(T)$, hence the
whole secondary fan, as lying in $L_\R^*$.

\subsec{toric ideal and the Gr\"obner fan}

We briefly review some results of Sturmfels{\sturmfels\sturmfelsII}.
Consider the generic orbit $T_M\cdot[1,..,1]$ in $\P\C^\cA$.
The closure of this orbit is a projective toric variety. Its
vanishing ideal $\cI_\cA$ is a homogeneous prime ideal called
a toric ideal, in $\C[y_0,..,y_p]$.
It is known that $\cI_\cA=\bra y^{l^+}-y^{l^-}|l\in L\ket$.
Here  $l^+_\mu=l_\mu$ if $l_\mu\geq0$ and zero otherwise, and 
similarly for $l^-$.

Every vector $\omega\in\R^{p+1}=\R^\cA$ induces a polyhedral subdivision 
$T_\omega$ of $\cA$ as follows. 
Consider the $(n+1)$-polytope $P_\omega:=conv(\{(\omega_\mu,
\mu)|\bar\mu\in\cA\})$
which is a lifting of $conv(\cA) = conv~(\{(1,\mu) | \bar\mu\in\cA\})$ 
according to the height vector $\omega$.
The lower envelope of $P_\omega$ is a polyhedral $n$-ball which maps
bijectively onto $conv(\cA)$. Let $T_\omega$ be the image under  this
projection. If $\omega$ is sufficiently generic $T_\omega$ is
a triangulation of $\cA$ (really of $conv(\cA)$). 
A triangulation $T$ of $\cA$ is regular
iff $T=T_\omega$ for some $\omega$. When $T$ is regular,
the set of all $\omega$ with $T=T_\omega$ is the 
interior $Int(\cC'(T))$ of the cone $\cC'(T)$ in the
secondary fan. 

Given a vector $\omega$, we define a weight on the
set of monomials in $y_0,..,y_p$, by
$t_w(y^l)=w_0l_0+\cdots+w_pl_p$. 
Consider now the leading term ideal $LT_\omega(\cI_\cA)$
of the toric ideal $\cI_\cA$, relative to the weight $t_\omega$.
 Two vectors $\omega,\omega'$
are said to be equivalent if $LT_\omega(\cI_\cA)=LT_{\omega'}(\cI_\cA)$. 
The equivalence classes of such vectors
form a complete fan in $\R^{p+1}$. Projecting
this fan along the linear subspace $\bar M_\R\hookrightarrow\R^{p+1}$
(given by $u\mapsto\sum\bra u,\bar\mu_i\ket e_i\in \R^{p+1}$), 
we get a fan in $L_\R^*=\R^{p+1}/\bar M_R$, which is called the
Gr\"obner fan $G\Sigma$ for the toric ideal
$\cI_\cA$. An interior point of a large cone
in $G\Sigma$ (or its lifting in $\R^{p+1}$)
is called a term order for $\cI_\cA$.

The Stanley-Reisner ideal $SR_T$ for a triangulation $T$ of $\cA$
is the ideal in $\C[y_0,..,y_p]$ 
generated by all monomials $y_I:=y_{i_1}\cdots y_{i_k}$ 
where $I=\{\bar\mu_{i_1},..,\bar\mu_{i_k}\}\subset\cA$
is not a simplex of $T$. 
A collection of vertices $\cP=\{\bar\mu_{i_1}, 
\cdots, \bar\mu_{i_k}\}$ is called primitive if $\cP$ is not 
a simplex contained in $T$ but any subset of it is a simplex. 
Then it's easy to see that $SR_{T}$ is generated
by the $y_\cP$ with $\cP$ being a primitive collection of
$T$.

\theorem{\sturmfels Let $\omega$ be a term order of the toric
ideal $\cI_\cA$. Then the polyhedral subdivision $T_\omega$ is a
triangulation of $\cA$. Moreover
$SR_{T_\omega}=rad~LT_\omega(\cI_\cA)$ and $G\Sigma$ is a refinement
of $S\Sigma$.}
\thmlab\SturmfelsThm

\subsec{toric variety $\P_\Sigma$ and Chow ring}

Let us consider a complete fan $\Sigma$ in $N$ 
with the property that the set of 
the primitive generators $\Sigma(1)$ is given by 
\eqn\dumb{\Sigma(1)=\{\mu \; | \; \bar\mu\in\cA, \, (\mu\not=0) \}.}
The fan $\Sigma$ need not be regular in general. 
In the following sections, we will restrict our attentions to the 
case in which $\Sigma$ is regular. In this case, the toric variety 
associated to the fan $\Sigma$ is a complete smooth toric 
variety $\P_\Sigma$. 

\theorem{(see {\oda}{\fulton}) 
For a complete smooth toric variety $\P_\Sigma$, 
the intersection 
ring $A^*(\P_\Sigma)$ is given by $\Z[D_1,\cdots,D_p]/{\cal I}$, 
where ${\cal I}$ is the ideal generated by 
\item{(i)} $D_{i_1}\cdots \cdot D_{i_k}$ for $\mu_{i_1},
\cdots,\mu_{i_k}$ not in a cone of $\Sigma$,  
\item{(ii)} $\sum_{i=1}^p \langle u, \mu_i \rangle D_i$ for $u\in M$.}

We call $\cI$ the Chow ideal of $\P_\Sigma$. Note that the similarity of 
(i) to the Stanley-Reisner ideal for a triangulation of $\cA$. Because of
this we call the ideal generated by (i) the Stanley-Reisner ideal of $\Sigma$. 
Equivalently,
this is the ideal generated by the elements $D_{i_1}\cdots D_{i_k}$ with 
$\{\mu_{i_1},..,\mu_{i_k} \}$ ranges over all primitive collections of the 
fan $\Sigma$ (, see below for the definition).
It follows from the theorem that $A^1(\P_\Sigma)=H^2(\P_\Sigma,\Z)
=\Z^p/M$.

\definition{(see \batyrevQ) A primitive collection of a complete 
fan $\Sigma$ is a subset $\cP\subset\Sigma(1)$ that 
does not generate a cone in $\Sigma$, but $\cP-\mu$ 
generates a cone in $\Sigma$ for every $\mu\in\cP$.}

Consider a piecewise linear function relative to $\Sigma$.
As before, the set of all piecewise linear funcions $PL(\Sigma)$ 
is canonically isomorphic to $\R^p=\R^{\Sigma(1)}$. The subset 
of the strictly convex piecewise linear functions will be 
denoted as $K(\P_\Sigma)$. Under our assumption of regularity 
of the fan $\Sigma$, we have 

\theorem{(see {\oda}) 
$PL(\Sigma)/M_\R\cong H^2(\P_\Sigma,\R)$, and
$K(\P_\Sigma)$ is the K\"ahler cone of $\P_\Sigma$. 
Moreover the first Chern class
$c_1(\P_\Sigma):=c_1(\Theta_{\P_\Sigma})$ is represented by 
$\alpha_\Sigma\in PL(\Sigma)$ with $\alpha(\mu)=1$ for all $\mu\in\Sigma(1)$.}

\theorem{(see {\batyrevQ}) 
$f\in K(\P_\Sigma)$ iff
$\sum_{\mu\in\cP}f(\mu)>f(\sum_{\mu\in\cP}\mu)$ for all primitive collections
$\cP$.}
\thmlab\odaparkthm

Let $\cP\subset\Sigma(1)$ be a primitive collection. Then
there is a unique cone $\sigma\in\Sigma$ of minimal dimension such that
$\sum_{\mu\in\cP}\mu\in Int(\sigma)$. By regularity $\sigma$ has a
unique set of integral generators $\cG\subset\Sigma(1)$. Thus we have a relation
\eqn\dumb{\sum_{\mu\in\cP}\mu=\sum_{\mu\in\cG}c_\mu\mu}
for some positive integral $c_\mu$ 
(positive because the LHS is in the interior; integral because
LHS is a lattice point).
Equivalently we can write
\eqn\prim{\sum_{\mu\in\cP}\bar\mu=c_0\bar\mu_0+\sum_{\mu\in\cG}c_\mu\bar\mu}
where $\bar\mu_0=1\times\vec0$ and $c_0:= |\cP|-\sum_{\mu\in\cG} c_\mu$.

\proposition{The sets $\cP$, $\cG\cup\{\vec0\}$ are disjoint. 
Hence \prim~defines an element $l(\cP)\in L$, 
which is called a primitive relation. If $l(\cP)_0\leq0$ then 
$l(\cP)^\pm$ are given by the left and right hand sides of 
{\prim}, respectively. }
\thmlab\PrimitiveRelationDefn
\proof First note that the fact that $\cP$ is primitive implies that
$\vec0\notin\cP$. Thus we need to show that $\cP,\cG$ are disjoint.
Suppose $\lambda$ is in both. Then we have
\eqn\stupid{
\sum_{\mu\in\cP-\lambda}\mu=\sum_{\mu\in\cG-\lambda}c_\mu\mu
+(c_\lambda-1)\lambda.}
But $\cP-\lambda$ generates a cone $\tau$ and so the LHS
of \stupid~ is in $Int(\tau)$. But we know the RHS is
in $\sigma$. This means that both sides are in $\tau\cap\sigma$.
So both sides are either in $\partial\tau$ or in $\partial\sigma$
or else $\tau=\sigma$.
It cannot be in $\partial\tau$ because it is in $Int(\tau)$.
We cannot have $\tau=\sigma$ because this would mean that 
$\cP-\lambda=\cG$ since these are their respective generating sets.
This would contradict $\lambda\in\cG$.

So finally suppose both sides are in $\partial\sigma$. 
This means that
$\tau\subset\partial\sigma$. This implies on the RHS of \stupid~
one of the coefficients of those generators must be zero. Now $c_\mu>0$
means that $c_\lambda-1=0$. Since the cones $\tau,\sigma$ are
regular, they have unique sets of generators. In particular
the generators of $\tau$ must be a subset of generators of $\sigma$,
ie. $\cP-\lambda\subset\cG$. But since both sides of \stupid~
have a unique nonnegative decomposition in terms of $\cG$, it follows
that $\cP-\lambda=\cG-\lambda$. Hence $\cP=\cG$ which contradicts
that $\cP$ is primitive. $\Box$

\proposition{$\bar K(\P_\Sigma)^\vee$ 
is generated by the primitive relations.}
\proof The natural pairing between $H_2(\P_\Sigma,\R)=L_\R$ 
and $H^2(\P_\Sigma,\R)=PL(\Sigma)/M_\R$ is given by
$L_\R\times PL(\Sigma)/M_\R\ra\R$,
$(l,f)\mapsto\sum_{\bar\mu\in\cA}l_\mu f(\mu)$.
Let $\cM\subset L_\R$ be the cone generated
by the primitive relations.
Given $f\in\bar K(\P_\Sigma)$ 
(modulo linear functions) and a primitive collection $\cP$, we have
$\sum_{\mu\in\cP} f(\mu)\geq f(\sum_{\mu\in\cP}\mu)$. 
So $\sum_{\mu\in\cP}f(\mu)-\sum_{\mu\in\cG}c_\mu f(\mu)\geq0$ or
$\sum l(\cP)_\mu f(\mu)\geq0$.
So every $l(\cP)$ is in
$\bar K(\P_\Sigma)^\vee$, ie. $\cM\subset\bar K(\P_\Sigma)^\vee$. We now show
that $\cM^\vee\subset\bar K(\P_\Sigma)$. 
Let $f\in\cM^\vee\subset PL(\Sigma)/M_\R$. Then for all primitive
collection $\cP$,
$\sum_{\bar\mu\in\cA}l(\cP)_\mu f(\mu)\geq0$, which
is equivalent, once again, to
$\sum_{\mu\in\cP} f(\mu)\geq f(\sum_{\mu\in\cP}\mu)$. 
Thus $f\in\bar K(\P_\Sigma)$.
$\Box$

\corollary{If $\tau$ is a cone contained in $\bar K(\P_\Sigma)$,
then $\tau^\vee$ contains all the primitive relations.}
\thmlab\TauPrimitiveRelation

We say that the fan $\Sigma$ has property (*) if it 
satisfies:
\eqn\dumb{\Sigma(1)\subset\partial conv(\Sigma(1))~~~~(*).}
This just means that the primitive generators of the one
dimensional cones in the fan $\Sigma$ lie on the faces of 
a convex polytope.
Note that each $(i+1)$-cone $\sigma\in\Sigma$ contains a canonical
$i$-simplex $s_\sigma$
whose vertices are the primitive generators of $\sigma$. 
Property (*) is also equivalent to the property that
each $s_\sigma$ lie on a boundary face of $conv(\Sigma(1))$. 

\theorem{The following are equivalent: 
\item{(i)} The first Chern class $c_1(\P_\Sigma)\in 
\bar K(\P_\Sigma)$.
\item{(ii)} Every primitive relation has $l(\cP)_0\leq0$.
\item{(iii)} The fan $\Sigma$ has property (*).
 }
\thmlab\chernclass
\proof We first show that (i) $\Leftrightarrow$ (ii).
Recall that $c_1(\P_\Sigma)$ is represented by the piecewise linear
function $f:=\alpha_\Sigma$ with $f(\mu)=1$ for all $\mu\in\Sigma(1)$.
By the above characterization of the K\"ahler cone, 
$c_1(\P_\Sigma)\in\bar K(\P_\Sigma)$ iff
for every primitive collection $\cP$ we have
\eqn\dumb{\sum_{\mu\in\cP}f(\mu)\geq f(\sum_{\mu\in\cP}\mu).}
The last condition is equivalent to $|\cP|\geq 
f(\sum_{\mu\in\cP}\mu)=f(\sum_{\mu\in\cG}c_\mu\mu)=\sum_{\mu\in\cG}c_\mu$. 
The last equality uses the linearity of $f$ in $cone(\cG)$. The 
inequality is in turn equivalent to 
$l(\cP)_0=\sum_{\mu\in\cG}c_\mu-|\cP|\leq0$.

We now show (i) $\Leftrightarrow$ (iii). Consider
the constant function which takes value 1 on $N_\R$. Its graph
in $\bar N_\R$ is the hyperplane $1\times N_\R$.
The fan $\Sigma$ which covers $N_\R$
can be translated vertically
to $\bar\Sigma$ which now covers $1\times N_\R$.
We denote by $s_{\bar\sigma}$ the translate of the simplex $s_\sigma$.
Now we pull the apex 
$\bar\mu_0=1\times\vec0$ down to the origin $0\times\vec0$ piecewise
linearly while holding each $s_{\bar\sigma}$ fixed.
The resulting piecewise linear hypersurface in $\bar N_\R$
is the graph of the
function $f$ representing $c_1(\P_\Sigma)$.
It's clear that $f$ is a convex function
iff $\bar\Sigma(1)$ lies on the boundary
of the convex polytope $conv~\bar\Sigma(1)$ in the hyperplane $1\times
N_\R$. Here $\bar\Sigma(1)$ is a translate of $\Sigma(1)$.
So $\bar\Sigma(1)$ lies on the boundary of $conv~\bar\Sigma(1)$
iff $\Sigma$ has property (*). 

This completes the proof.
$\Box$

\newsec{Compactification and Existence of Maximal Degeneracy Points}

Now we consider a compactification of $Hom(L,\C^\times)$ in
proposition {\PrincipalBundle} 
using the two complete fans $S\Sigma$ and $G\Sigma$ in 
$L_\R^*=Hom(L,\Z^{p-n})\otimes\R$. Under our assumption that $\Sigma$ is a 
regular fan, we have the following exact sequence 
\eqn\dumb{0\ra L\ra\Z^\cA\ra\bar N\ra0}
where $\Z^\cA\ra\bar N$ is given by 
$l\mapsto \sum_{\bar\mu\in\cA} l_\mu\bar\mu$. This
is onto because $\Sigma(1)$ generates the lattice $N$
since $\Sigma$ is regular. 
Dualizing the map $\Z^{\cA}\ra \bar N$, we get $\bar M\ra\Z^{\cA}$,
$\bar u\mapsto\sum\bra \bar u,\bar\mu\ket e_{\bar\mu}$ where 
the $e_{\bar\mu}$
is the standard basis of $\Z^{\cA}$. 
Then we have the isomorphism 
\eqn\latiso{ 
L^* \simeq \Z^{\cA}/{\bar M} \;\;. }
With this integral structure, both $S\Sigma$ and $G\Sigma$ in $L^*$ 
consist of rational polyhedral cones and thus we can consider the toric
varieties $\P_{S\Sigma}$ and $\P_{G\Sigma}$. According to
Theorem {\SturmfelsThm}, 
$\P_{G\Sigma}$ is a partial resolution of $\P_{S\Sigma}$ which is typically
singular.

\problem{(Existence) Given an equivariant resolution $\phi:\P_{G\Sigma'}\ra
\P_{G\Sigma}$, when does $\P_{G\Sigma'}$ admit a maximal degeneracy point?}

We will prove the following:  

\theorem{If the regular fan $\Sigma$ has the property (*) 
and the (K\"ahler) cone $K(\P_\Sigma)$ is a large cone in 
$L^*_\R$ then for every equivariant resolution 
$\phi:\P_{G\Sigma'}\ra \P_{G\Sigma}$,
there is at least one maximal degeneracy point in
$\P_{G\Sigma'}$.}
\thmlab\MainTheorem

\remark{As we have proved, the combinatorial conditions
in the hypotheses of the theorem are equivalent to that
the toric variety $\P_\Sigma$ is smooth, projective,
and has $c_1(\P_\Sigma)\in\bar K(\P_\Sigma)$.}

\subsec{indicial ideals}

Let $\tau$ be a regular large cone in the space $L_\R^*$.
Then it determines a unique
$\Z$-basis $\{ l^{(1)},\cdots,l^{(p-n)} \}$
of $L$ in $\tau^\vee \cap L$, and hence a set of canonical coordinates
$x^{(k)}_\tau$ of the affine variety 
$U_\tau=Hom(\tau^\vee\cap L,\C)$. 
Thus by Proposition \PrincipalBundle, we have

\proposition{For any regular large cone $\tau$ in $L_\R^*$,
we have an embedding $(\C^\times)^\cA/(\C^\times\times T_M)\ra
Hom(\tau^\vee\cap L,\C)$, given by $a\mapsto x_\tau$ where
$x^{(k)}_\tau=(-1)^{l_0^{(k)}} a^{l^{(k)}}$. }
\thmlab\UpstairDownstair

Recall that the interior of a large cone $\tau$ in the Gr\"obner fan $G\Sigma$
is an equivalence class of
 term order $\omega\in\R^{p+1}=\R^\cA$ for the toric ideal $\cI_\cA$.
We identify $\cI_\A$ in $\C[y_0,\cdots,y_p]$ 
with the ideal generated by $\{ \Box_l \}_{l\in L}$ in 
$\IC[{\pd \; \over \pd a_0},\cdots, {\pd\; \over \pd a_p} ]$.
Let $B_\omega$ be a minimal Gr\"obner basis
(see {\sturmfelsII}) for the term order $\omega$.
Then $B_\omega$ is
a finite set of $\Box_l$ with $\omega\cdot l\geq0$ for all $\omega\in\tau$.
Moreover the leading term relative to $\omega$
is $LT_\omega\Box_l=({\pd \; \over \pd a})^{l^+}$. Now
\eqn\msqbox{
a^{l^+} \Box_l = a^{l^+}\( {\pd \; \over \pd a} \)^{l^+} - 
a^{l^+ - l^-} a^{l^-} \( {\pd \; \over \pd a} \)^{l^-}.
}
Since $\omega\cdot(l^+-l^-)\geq0$ for all $\omega\in\tau$, we have 
$l=l^+-l^- \in \tau^\vee\cap L$.
Since $\tau^\vee\cap L$ is generated by 
$\{l^{(1)},..,l^{(p-n)}\}$, 
it follows form Proposition \UpstairDownstair~ that  
$a^{l^+ -l^-}$ in the second term  
can be expressed as a monomial of $\{ x_\tau^{(k)} \}$ which vanishes when 
$x_\tau^{(k)} \rightarrow 0$. The first term in \msqbox\ are 
`homogeneous' and can be written as a constant coefficient
 polynomial $I_l(\ta{0},\cdots,\ta{p})$ in the log derivatives
$\theta_{a_i}=a_i {\pd \; \over \pd a_i }$. 
Thus the leading term of $a^{-\gamma} a^{l^+}\Box_la^\gamma(1+O(x))$
is obviously $I_l(\gamma)$.

With the above motivation we define for any $l\in L$
\eqn\dumb{I_l(\gamma):=a^{-\gamma} a^{l^+}
({\partial\over\partial a})^{l^+}a^\gamma\in
\C[\gamma]:=\C[\gamma_0,..,\gamma_p].}

\definition{For a cone $\tau\subset L_\R^*$ and
an exponent $\beta=-1\times\vec0 \in\bar N_\R$, 
the indicial ideal $Ind(\tau,\beta)$
is the ideal in $\C[\gamma]$ generated by
\eqn\dumb{I_l(\gamma)~with~l\in\tau^\vee\cap L,~ l\neq0,~~~~
\sum\bra u,\bar\mu\ket\gamma_\mu-\bra u,\beta\ket~with~u\in \bar M.} }
\thmlab\IndicialIdeal

Note that the second set of relations is equivalent to
$\sum_{\bar\mu\in\cA}\bar\mu\gamma_\mu-\beta$. Note also that we
make no assumption on the cone $\tau$.
In practice, one can sometimes reduce the problem of finding
the indices (ie. the zero locus of the indicial ideal) by using
a Gr\"obner basis of the toric ideal. 

\proposition{Let $\tau$ be a large cone in the Gr\"obner fan 
and $B=\{\Box_{l^{(1)}},..,
\Box_{l^{(g)}}\}$ be a minimal Gr\"obner basis of $\cI_\cA$ for
the class of term orders $\omega\in\tau$. Then the ideal
\eqn\dumb{\cJ_B:=\bra I_{l^{(i)}}(\gamma),~i=1,..,g; ~\sum\bra
u,\bar\mu\ket\gamma_\mu-\bra u,\beta\ket~\ket}
has the same zero locus as does $Ind(\tau,\beta)$.}
\proof By definition for each $i$, we have
$LT_\omega\Box_{l^{(i)}}=({\partial\over\partial a})^{l^{(i)+}}$
and $\bra\omega,l^{(i)}\ket>0$ for a term order $\omega\in\tau$.
Thus $l^{(i)}\in\tau^\vee \cap L$. In particular $\cJ_B\subset
Ind(\tau,\beta)$.

Suppose $\cJ_B$ vanishes at $t=(t_0,..,t_p)\in\C^{p+1}$. We must show
that so does $Ind(\tau,\beta)$. It is enough to check that
$I_l(t)=0$ $\forall l\in\tau^\vee\cap L\backslash(0)$. 
Let us fix $l$. 
Then any interior $\omega\in\tau$ has $\bra\omega,l\ket>0$.
Since $B$ is a Gr\"obner basis, we have
$LT_\omega\Box_l=\sum b_i LT_\omega\Box_{l^{(i)}}$ for some
differential operators $b_i$. This gives
$({\partial\over\partial a})^{l^+} a^t=\sum b_i 
({\partial\over\partial a})^{l^{(i)+}} a^t=0$.
This implies that $I_l(t)=0$.
 $\Box$

\subsec{finite codimension theorem}

\theorem{Suppose the regular fan $\Sigma$ has the property (*). 
If $\tau$ is any cone contained in $\bar K(\P_\Sigma)$,
then there is a canonical onto ring homomorphism
$A(\P_\Sigma)\otimes\C\ra\C[\gamma]/Ind(\tau,\beta)$
with $D_i\mapsto\gamma_i$.}
\thmlab\OntoRingMap
\proof 
By Theorem \chernclass, we have $l(\cP)_0\leq0$ for each primitive
collection $\cP$. By Proposition \PrimitiveRelationDefn, 
it follows that $l(\cP)^+_0=0$ and that
$l(\cP)^+_i=1$ if $\mu_i\in\cP$ and 0 otherwise. Thus
we have
 $I_{l(\cP)}(\gamma)=a^{-\gamma} a^{l(\cP)^+}
({\partial~\over\partial a})^{l(\cP)^+} a^\gamma=\gamma^{l(\cP)^+}$.
As remarked earlier the elements
 $I_{l(\cP)}(D)=D^{l(\cP)^+}$, as $\cP$ ranges over the primitive
collections, generate
precisely the Stanley-Reisner ideal of the fan $\Sigma$.

Consider the isomorphism
$\C[\gamma_0,..,\gamma_p]/\bra\sum\gamma_i+1\ket\ra\C[D_1,..,D_p]$,
$\gamma_i\mapsto D_i$ for $i\neq0$, $\gamma_0\mapsto-1-D_1-\cdots-D_p$.
For $\beta=-1\times\vec0$, 
it is clear that under this isomorphism we have
\eqn\dumb{\eqalign{
\bra\sum\bra u,\bar\mu_i\ket\gamma_i-\bra u,\beta\ket,~u\in \bar M\ket=&
 \bra\sum\gamma_i+1,~\sum\bra u,\mu_i\ket\gamma_i,~u\in M\ket\cr
&\ra\bra\sum\bra u,\mu_i\ket D_i,~u\in M\ket.}}
Comparing with the generators (i) of the Chow ideal, we get an isomorphism
\eqn\dumb{\C[\gamma]/\bra I_{l(\cP)}(\gamma),~\cP~primitive;~
\sum\bra u,\bar\mu_i\ket\gamma_i-\bra u,\beta\ket,~u\in\bar M\ket
\ra A(\P_\Sigma)\otimes\C.}

By Definition \IndicialIdeal~ and 
Corollary \TauPrimitiveRelation, $Ind(\tau,\beta)$ contains $I_{l(\cP)}(\gamma)$
for all primitive collection $\cP$, and all
$\sum\bra u,\bar\mu_i\ket\gamma_i-\bra u,\beta\ket,~u\in M$.
Thus inverting the above isomorphism, we get an onto homomorphism
$A(\P_\Sigma)\otimes\C\ra\C[\gamma]/Ind(\tau,\beta)$. $\Box$

\corollary{Under the hypotheses of Theorem \OntoRingMap,
$Ind(\tau,\beta)$ has finite codimension in $\C[\gamma]$.}
\proof This follows from the fact that $A(\P_\Sigma)\otimes\C$ is finite
dimensional. $\Box$

\corollary{Under the hypotheses of Theorem \OntoRingMap, 
the zero locus of  the 
ideal $\bra I_{l(\cP)}(\gamma),~\cP~primitive;~
\sum\bra u,\bar\mu_i\ket\gamma_i-\bra u,\beta\ket,~u\in\bar M\ket$ 
has exactly one point $\gamma=(-1,0,..,0)\in\R^{p+1}$.}
\proof  As pointed out in the proof of the theorem,
the relations in this ideal are
\eqn\IdealRelation{\eqalign{
I_{l(\cP)}(\gamma)=&0~~with~\cP~primitive\cr
\sum\bra u,\mu_i\ket\gamma_i=&0~~with~u\in M\cr
\sum\gamma_i+1=&0.}}
The first two sets of relations are identical to those of the Chow ideal.
Since the Chow ideal is homogeneous and has finite codimension, the first
two sets of relations above give $\gamma_1=\cdots=\gamma_p=0$.
It follows from the third relation that $\gamma_0=-1$. $\Box$

\corollary{Under the hypotheses of Theorem \OntoRingMap,
the zero locus of the ideal $Ind(\tau,\beta)$ has at most the one
point $\gamma=(-1,0,..,0)$.}
\proof This follows from Theorem \OntoRingMap~ and the preceding corollary.
$\Box$

\lemma{Let $S$ be a nonempty set of integral points in a large regular cone
$C$ in $L_\R^*$. Then there exists a strongly convex regular cone $C'$
containing both $C$ and the set $S-\delta$ for some $\delta\in S$.}
\proof Without loss of generality, we identify $C$ with the
cone $\R^{p-n}_{\geq 0}$. It is enough to find a
strongly convex cone $C'$ (not necessary regular)
with the above inclusion property, for every such cone
is contained in a strongly convex regular cone. 
Now if $S$ contains the origin, then we set $C'=C$.
Otherwise we proceed as follows.
Given $x\in C$, $v\in Int(C^\vee)$,
let $H_v(x)$ be the hyperplane in $\R^{p-n}$ containing
$x$ and being normal to $v$. We fix $x,v$ and 
let $t$ be the minimal positive number such that
$H_v(tx)$ meets $S$. 
By perturbing $v$, if necessary, we may assume
that $H_v(tx)$ meets $S$ at exactly one point $\delta$.
In particular we have $\bra v,s-\delta\ket>0$ for
all $s\in S\backslash \{\delta\}$. Performing a parallel translation of 
the hyperplane to $H_v(0)$ we see that the points of both $C$ and $S-\delta$,
except the origin, all lie in the same open half space  bounded by
$H_v(0)$. It follows that the normal cone to $C\cup(S-\delta)$
at $0$ is strongly convex, and obviously has the desired
inclusion property. $\Box$.

\corollary{Under the hypotheses of Theorem \OntoRingMap,
if $\tau$ is a regular large cone in $L_\R^*$ then the GKZ system
has at most one power series solution of the form
$a^\gamma(1+g(x))$ where $g(0)=0$. (The $x$ are
the canonical coordinates of the 
affine variety $U_\tau=Hom(\tau^\vee\cap L,\C)$.)
Moreover if this is
a solution then $\gamma=(-1,0,..,0)$.}
\thmlab\UniquenessIndex
\proof The second statement follows from the preceeding corollary.
For the first statement we show that any solution of the form
$a^\gamma f(x)$ with $f(0)=0$ is identically zero.
Suppose not. Then we can write $f(x)=\sum_{l\in\tau^\vee} \lambda_l a^l$
with $\lambda_l$ nonzero for some $l$.
Let $S:=\{l\in\tau^\vee|\lambda_l\neq0\}$. Since $f$
has no constant term, $S$ does not contain the origin.
By the preceeding lemma, we have a strongly convex regular
cone $C'$ containing both $\tau^\vee$ and $S-\delta$ for some
$\delta\in S$. Now put $\lambda_l \equiv 0$ for $l$ outside 
$\tau^\vee$. We have
\eqn\dumb{
a^\gamma f(x)=a^\gamma\sum_{l\in S}\lambda_l a^l
=a^{\gamma+\delta}\sum_{l\in C'}\lambda_{l+\delta} a^l.}
But the last sum can be written as $\lambda_\delta a^{\gamma+\delta}(1+g(x'))$
where $g(0)=0$, and the $x'$ are the canonical coordinates
of $U_\sigma (\sigma={C'}^\vee)$. 
By the preceeding corollary applied to $\sigma\in\tau$, 
we have $\gamma+\delta=\gamma$, which means that $0=\delta\in S$.
This contradicts the fact that $S$ does not contain the origin. $\Box$

The ideal $Ind(\tau,\beta)$ can be all of $\C[\gamma]$ in
which case its zero locus is empty. For example
if $\tau=(0)$ and $l_1\mu_1+\cdots+l_p\mu_p=0$ for some nonpositive integers
$l_1,..,l_p$, not all zero, then $l=(l_0:=-l_1-\cdots-l_p,l_1,..,l_p)\in L$
and we have $I_l(-1,0,..,0)\neq0$.

\corollary{Under the hypotheses of Theorem \OntoRingMap,
the zero locus of the ideal $Ind(\tau,\beta)$ has exactly the point
$\gamma=(-1,0,..,0)$ iff every nonzero
$l\in\tau^\vee\cap L$ has $l_i>0$ for some $i>0$.}
\proof It's obvious that the point $\gamma=(-1,0,..,0)$
 satisfies the second and third relations
in \IdealRelation. So we need only to consider the conditions
$I_l(-1,0,..,0)=0$ for all $l\in\tau^\vee\cap L$, $l\neq0$.
By definition, $I_l(-1,0,..,0)=a_0 a^{l^+}({\partial~\over\partial a})^{l^+}a_0^{-1}$. 
It is obvious that this is zero iff $l_i>0$ for some $i>0$. $\Box$

\subsec{proof of the main theorem \MainTheorem}

In this section, we will assume that 
the fan $\Sigma$ is regular and has property (*), and that
$\bar K(\P_\Sigma)$ is a large cone.
It is obvious that the fan endows
the polytope $conv(\cA)$ with a natural triangulation $T^0$
whose $n$-simplices correspond 1-1 to the large cones in $\Sigma$. 
We will call this a maximal triangulation $T^0$. 
Then  under the identifications 
$\R^{p+1}/\bar M_\R=\R^p/M_\R=L_\R^*$,
we have $\cC(T^0)=\bar K(\P_\Sigma)$.
Since $\bar K(\P_\Sigma)$ is
a large cone, the triangulation $T^0$ is regular. 

We can regard the cone $\bar K(\P_\Sigma)$ together with all
its faces as a fan, which we denote $K$. 
Corresponding to this fan is an affine toric variety
\eqn\dumb{\P_K=Hom(\bar K(\P_\Sigma)^\vee\cap L,\C).}
It is of course canonically a toric subvariety of $\P_{S\Sigma}$.
In $\P_K$ lies a unique torus fixed point $p_K$. In general $p_K$
is singular. Let $\P_{K'}\ra\P_K$ be any equivariant
resolution (thus $K'$ is a regular refinement of $K$).
The fiber over $p_K$ is a union of torus orbits.
{}From standard toric geometry, we get

\lemma{Given $\P_{K'}\ra\P_K$, the
torus fixed points in the fiber over $p_K$ correspond
1-1 with the large cones $\tau\in K'$.}

We denote by $p_\tau$ the fixed point corresponding to $\tau$.
As pointed out earlier, the period integrals $\bar\Pi_\alpha$
are local sections of a sheaf over the open variety $S_\cA$;
and the latter is canonically an open subvariety
in every toric variety $\P_F$ corresponding
to any rational fan $F$ defined in $L_\R^*$.
In particular for the fan $K'$,
we can ask for the behavior of the period integrals near $p_\tau\in
\P_{K'}-S_\cA$.

Recall that we have an explicit formula for a period integral
$\bar\Pi_{\alpha_0}$ given by Proposition \PowerSeries. 
Near $p_\tau$, we can write this period in the canonical
coordinates $x=x_\tau$ using Proposition \UpstairDownstair.

\proposition{Fix a large cone $\tau\in K'$ and let
$l^{(1)},..,l^{(p-n)}\in L$ be the canonical generators of $\tau^\vee$.
Near $p_\tau$ in the canonical coordinates $x=x_\tau$,
the period integral
$\bar\Pi_{\alpha_0}$ can be written as
\eqn\dumb{
\bar\Pi_{\alpha_0}=\sum_{m_1,\cdots,m_{p-n} \geq 0,~\sum m_k l^{(k)}_0\leq0}
  {\Gamma(-\sum m_k l^{(k)}_0 +1) \over
     \prod_{1\leq i \leq p} \Gamma(\sum m_k l^{(k)}_i +1) }
	x^m. }
In particular it extends holomorphically across $p_\tau$.} 
\proof 
Since $\tau\subset\bar K(\P_\Sigma)=\cC(T^0)$,
we have $\tau^\vee\supset\cC(T^0)^\vee=\bar\cC'(T^0)^\vee$. 
Thus it follows from
Proposition \GKZprop~ that $\tau^\vee$ contains
the $K(I):=\{x\in L_\R|x_i\geq0,~i\notin
I\}$ for all $I\in T^0$. Since $T^0$ is the triangulation of $conv(\cA)$
induced by the fan $\Sigma$, every $n$-simplex of $T^0$
has the point $\bar\mu_0=1\times\vec0$ as a vertex, ie.
every basis $I\in T^0$ contains this point.
It follows that for every tuple $(l_1,..,l_p)$
we sum over in Proposition \PowerSeries,
if we set $l_0=-l_1-\cdots-l_p$ then we get
$l:=(l_0,..,l_p)\in K(I)\subset \tau^\vee$.
So there exists unique integers $m_1,..,m_{p-n}\geq0$ such that
$l_i=\sum m_k l^{(k)}_i$ for all $i$. Thus
$(-1)^{l_1+\cdots+l_p}
{(l_1+\cdots+l_p)!\over l_1!\cdots l_p!}
{a_1^{l_1}\cdots a_p^{l_p}\over a_0^{l_1+\cdots+l_p}}$ coincides
with ${\Gamma(-\sum m_k l^{(k)}_0 +1) \over
   \prod_{1\leq i \leq p} \Gamma(\sum m_k l^{(k)}_i +1) } x^m$.
   Conversely given a nonnegative tuple $(m_1,..,m_{p-n})$ with
   $\sum m_k l^{(k)}_0\leq0$, the last term is either zero
   (when some $\Gamma(\sum m_k l^{(k)}_i +1)$ is infinite), or
   it coincides with some
$(-1)^{l_1+\cdots+l_p}
{(l_1+\cdots+l_p)!\over l_1!\cdots l_p!}
   {a_1^{l_1}\cdots a_p^{l_p}\over a_0^{l_1+\cdots+l_p}}$ when nonzero.
   $\Box$

Thus on the one hand we have one period integral which
extends across $p_\tau$. On the other hand there can't be
any other period integrals with this property
by Corollary \UniquenessIndex. So we conclude that

\theorem{For any equivariant resolution $\P_{K'}\ra\P_K$,
near a given torus fixed point in $\P_{K'}$, 
$\bar\Pi_{\alpha_0}$ extends holomorphically
across that point and is the only period
integral which does so.}

\corollary{For any equivariant resolution 
$\phi:\P_{S\Sigma'}\ra\P_{S\Sigma}$ 
near a given torus fixed point in $\phi^{-1}(p_K)$, 
$\bar\Pi_{\alpha_0}$ extends holomorphically
across that point and is the only period
integral which does so.}
\proof Since $\bar K(\P_\Sigma)$ is a large cone in the secondary
fan $S\Sigma$, we can regard the fan $K$ (consisting of all the
faces of $\bar K(\P_\Sigma)$) as a subfan of $S\Sigma$.
The preimage $K':=\phi^{-1}(K)\subset S\Sigma'$ 
is a regular fan refinement of $K$, and
the induced map $\P_{K'}\ra\P_K$ is a restriction of 
$\phi:\P_{S\Sigma'}\ra\P_{S\Sigma}$. Applying the preceding theorem,
we get the desired result. $\Box$

Since the Gr\"obner fan $G\Sigma$ is a refinement of the
secondary fan $S\Sigma$, any equivariant resolution
$\P_{G\Sigma'}\ra\P_{G\Sigma}$ is also a resolution of
$\P_{S\Sigma}$. Thus the preceding corollary applies.
This proves our main theorem \MainTheorem.

\remark{Actually we have proved a stronger result than \MainTheorem.
Not only is $\bar\Pi_{\alpha_0}$ the only period integral that extends
across a torus fixed point given above, this is the only (up to
multiple) solution to the GKZ system which extends.}

\newsec{Gr\"obner bases}

As we have seen, we can obtain rather detailed information about
how period integrals, or more generally the solutions
to the GKZ system behave near a particular boundary point
by studying a cone in the secondary fan. In the preceding discussion
it was the K\"ahler cone and the corresponding torus
fixed point  which play an essential role.

There is in fact more to the story. We can actually reduce the GKZ system
to a finite systems of PDEs in a {\it canonical} way. These finite
systems correspond 1-1 to the Gr\"obner bases of the toric ideal.
The Gr\"obner fan $G\Sigma$ turns out to be just the right set up for
constructing this finite systems.
One can then use them to study explicit solutions
to the GKZ system near various boundary component in the parameter space.
Near a maximal degeneracy point given above, the finite system 
turns out to have a uniform and simple description in terms 
of primitive collections of
the fan $\Sigma$.

For generic values of the exponent $\beta$
(the so-called nonresonant cases), 
the theory of Gel'fand-Kapranov-Zelevinsky
\GKZ~ gives a uniform description for the solutions to
the GKZ system. However for our special values of the exponent 
$\beta=-1\times 0$ an analogous uniform description is not known.
This exceptional case is especially important for applications in
mirror symmetry. In this case, there are conjectural formulas for
general solutions to the GKZ system{\HKTYII}{\HLY}. Using our construction
of finite PDE system, we can explicitly verify these solution
formulas.

\subsec{the maximal triangulation $T^0$}

As in the previous section, we assume the fan $\Sigma$ is regular 
and has the property (*). We also assume the cone $K(\P_\Sigma)$ 
is large. 
As we have seen at the beginning of 
the section 5.2, we can associate the fan $\Sigma$ to the 
maximal triangulation $T^0$ of $\cA$. Then the cone 
$\cC'(T^0)$ in the secondary fan coincides with
the closure  $\bar K(\P_\Sigma)$

\proposition{For any term order $\omega$ with $T_\omega=T^0$, we have
$SR_{T^0}=LT_\omega(\cI_\cA)$, ie. $LT_\omega(\cI_\cA)$ is radical.}
\proof By Theorem \chernclass, $l(\cP)_0\leq0$ for every
primitive collection $\cP$ of the fan $\Sigma$.
Since $T_\omega=T^0$, we have $\omega\in Int~\cC'(T^0)$. Thus $\omega$ modulo
$\bar M_\R$ is an element $f$ in the K\"ahler cone 
$K(\P_\Sigma)\subset PL(\Sigma)/M_\R$. This is the piecewise linear
function given by
$f(\mu_i)=\omega_i-\omega_0$.
By Theorem \odaparkthm, for every primitive collection $\cP$
of $\Sigma$ we have
\eqn\fcond{\sum_{\mu\in\cP}f(\mu)>f(\sum_{\mu\in\cP}\mu).}
If we write $\sum_{\mu\in\cP}\mu=\sum_{\mu\in\cG}c_\mu\mu$ where $\cG$
is the generator set of the smallest cone in $\Sigma$ containing
$\sum_{\mu\in\cP}\mu$, then \fcond~ reads $\bra\omega,l(\cP)\ket>0$.
This implies that $LT_\omega(y^{l(\cP)^+}-y^{l(\cP)^-})=y^{l(\cP)^+}$.
Since $l(\cP)_0\leq0$, it follows that $l(\cP)^+_i=1$ if $\mu_i\in\cP$
and 0 otherwise. Hence $y^{l(\cP)^+}=y_\cP$.
Since $SR_{T^0}=\bra y_\cP|\cP~primitive\ket$, it follows that
$SR_{T^0}\subset LT_\omega(\cI_\cA)$. By Theorem \SturmfelsThm,
we have the reverse inclusion. $\Box$

\corollary{The cone $\bar K(\P_\Sigma)$ is in the Gr\"obner fan
$G\Sigma$.}
\thmlab\UnrefinedThm
\proof We must show that any two term orders $\omega,\omega'$ in
$K(\P_\Sigma)=Int~\cC'(T^0)$ are equivalent. 
Now they being in $\cC'(T^0)$ means that $T_\omega=T^0=T_{\omega'}$. 
It follows from the preceding proposition that 
$LT_\omega(\cI_\cA)=LT_{\omega'}(\cI_\cA)$.
$\Box$

\corollary{The set $\{y^{l(\cP)^+}-y^{l(\cP)^-}|\cP~primitive\}$ is
a minimal Gr\"obner basis of the toric ideal $\cI_\cA$
for a term order $\omega$ with $T_\omega=T^0$.}
\thmlab\GrobnerBasis
\proof From the proof of the above proposition, we see that
the leading terms $y_\cP$
{}from the given set generate $LT_\omega(\cI_\cA)$,
and so this set is a Gr\"obner basis. Moreover the leading terms
have coefficient 1; and $y_\cP$ from one primitive collection
$\cP$ does not divide $y_{\cP'}$ from another $\cP'$ unless 
$\cP=\cP'$. $\Box$

According to Corollary \UnrefinedThm,
the preimage of the torus fixed point $p_K\in\P_{S\Sigma}$
under the blow-up $\P_{G\Sigma}\ra\P_{S\Sigma}$ is a single point. We 
denote it as $p_K$ also.

\corollary{Let $\phi:\P_{G\Sigma'}\ra\P_{G\Sigma}$ be an equivariant
resolution of the 'Gr\"obner variety\'. Every torus fixed point
in $\phi^{-1}(p_K)$ is a maximal degeneracy point.}

\conjecture{Given 
$\phi:\P_{G\Sigma'}\ra\P_{G\Sigma}$ as above, the only
maximal degeneracy points
should be torus fixed points in
$\phi^{-1}(p_K)$.}

\subsec{the GKZ system and primitive collections}

If we identify $\C[y_0,..,y_p]$ with the ring of differential
operators $\C[{\partial\over\partial a_0},..,{\partial\over\partial a_p}]$
in an obvious way, then the toric ideal $\cI_\cA$ becomes precisely the
the ideal generated by the order $>1$ operators in the GKZ system.
Thus a Gr\"obner basis {\CoxLittleOshea} for the toric ideal allows us to reduce the GKZ system to a finite system in a canonical fashion.

For example,
fix an equivariant resolution $\phi:\P_{G\Sigma'}\ra\P_{G\Sigma}$.
We consider the GKZ system near the torus fixed point 
$p_\tau\in\phi^{-1}(p_K)$.
Relative to the term orders $\omega$ with $T_\omega=T^0$,
Corollary \GrobnerBasis~ reduces the GKZ system to 
the finite system of equations
\eqn\dumb{a^{-\gamma}\Box_{l(\cP)}a^{\gamma}f(x)=0}
indexed by the primitive collections $\cP$ of $\Sigma$.
Here $\gamma$ satisfies $\sum\gamma_\mu\bar\mu=\beta$,
and $x=x_\tau$ is the canonical coordinates of the
affine subvariety $U_\tau\subset\P_{G\Sigma'}$.

Near other torus fixed points, not necessarily in $\phi^{-1}(p_K)$,
the GKZ system can also be reduced to a finite system using Gr\"obner
bases of the toric ideal
relative to other term orders. There is a degree bound
for the binomials appearing in such Gr\"obner bases \sturmfels. There is also
a combinatorial description for the universal Gr\"obner basis.
We have worked out some
examples of this in {\HLY}.

\subsec{mirror symmetry}

As we pointed earlier,
the general theory we have discussed arises naturally
in mirror symmetry. Let's first consider the famous example 
studied by
Greene-Plesser {\GreenePlesser} and Candelas et al 
\candelasETAL. Let $\nabla$ be the
standard $4$-simplex in $\R^4$ with vertices at the basis point
$e_1,..,e_4$ and at $-e_1-\cdots-e_4$. Let $\Sigma$ be the
fan consisting of cones over the faces of $\nabla$.
Then the toric variety $\P_\Sigma$ coincides with the projective
space $\P^4$. 
The set $\cA:=\{1\times\mu|\mu\in\nabla\cap\Z^4\}$ has 6 points.
Let $\Delta$ be the polar dual of $\nabla$.
The fan over the faces of $\Delta$ is not regular but admits a 
regular refinement, which we denote $\Phi$.
Let $\cL=K^{-1}_{\P_\Phi}$ be the anticanonical bundle of $\P_\Phi$.
The space of sections $H^0(\cL)$ is isomorphic to $\C^\cA$,
and the 5-dimensional projective space $\P H^0(\cL)=\P\C^\cA$
parameterizes a family of Calabi-Yau hypersurfaces $Z_a$
in $\P_\Phi$. Here $Z_a$, 
$a\in\P\C^\cA$, is the hypersurface given by
the closure in $\P_\Phi$ of the zero locus of $f_\cA(X,a)$.
Then via the Poincar\'e residue map, the periods of the $Z_a$ 
can be written as 
period integrals
\eqn\dumb{a_0\int_\alpha{1\over f_\cA(X,a)}\prod{dX_i\over X_i}.}
The Gr\"obner variety $\P_{G\Sigma}$ in this
case is $\P^1$. It's easy to show that with a suitable choice of
coordinate, $x=0$ is the unique maximal degeneracy point in $\P^1$.
A remarkable discovery in \candelasETAL~ is that a particular
combination of the period integrals produces a generating function
which predicts the number of rational curves in each degree in
a generic quintic $Y$ in $\P_\Sigma=\P^4$. A quintic $Y$ in $\P^4$
and a Calabi-Yau hypersurface $Z$ in $\P_\Phi$ make what's known
as a {\it mirror pair}, for the reason that their Hodge diamonds are
mirror reflections of one another. 

There is a straightforward generalization of the above example to $n\geq4$.
More generally Batyrev have given a construction of 
pairs of Calabi-Yau hypersurfaces from pairs of
reflexive polytopes generalizing the pair $\nabla,\Delta$ above
and all previously known examples. An important  question arising
{}from each such pair is the 
{\it question of existence of maximal degeneracy
points} in the deformation space of the hypersurfaces.
Our general theory for studying period integrals specializes to
these cases, 
settles the question of existence of maximal degeneracy points,
 and has given a uniform
construction of these points.
The existence and construction  of such points are crucial
for the computation of (the analogues) of the Candelas et al generating
function. We have verified in {\HLY} that our construction in
all previously known plus some new examples, correctly
predicts this generating function.

\vfill\eject   

\centerline{\bf References }

\item{\batyrevJA} \refbatyrevJA
\item{\batyrevDuke} \refbatyrevDuke
\item\batyrevQ \refbatyrevQ
\item{\BKK} \refBKK 
\item{\BFS} \refBFS 
\item{\candelasETAL} \refcandelasETAL
\item{\candelasIIa} \refcandelasIIa
\item{\candelasIIb} \refcandelasIIb
\item{\CoxLittleOshea} \refCoxLittleOshea
\item{\Ft} \refFt
\item{\fulton} \reffulton
\item{\GKZ} \refGKZ
\item{\GKZII} \refGKZII
\item{\GreenePlesser} \refGreenePlesser 
\item{\HKTYI} \refHKTYI
\item{\HKTYII} \refHKTYII
\item{\HLY} \refHLY
\item{\KT} \refKT
\item{\morrison} \refmorrison
\item{\MFK} \refMFK 
\item{\oda} \refoda
\item{\OdaPark} \refOdaPark 
\item{\sturmfels} \refsturmfels
\item{\sturmfelsII} \refsturmfelsII
\item{\yau} \refyau


\bye